\documentclass[12pt,oneside,reqno]{amsart}
\usepackage{a4}
\usepackage{amssymb}

\usepackage{verbatim}

\usepackage{lineno}

\numberwithin{equation}{section}

\def\wid{\mathrm{wid}}
\def\cF{\mathcal F}
\def\cB{\mathcal B}
\def\cC{\mathcal C}
\def\cM{\mathcal M}
\def\cX{\mathcal X}
\def\cY{\mathcal Y}

\def\l{\lambda}
\def\vk{\varkappa}
\def\f{\varphi}
\def\s{\sigma}

\def\N{\mathbf N}
\def\Z{\mathbf Z}
\def\R{\mathbf R}

\def\supp{\mathrm{supp}}
\def\str#1{\langle#1\rangle}
\def\sym#1{\mathrm{Sym}(#1)}
\def\aut#1{\mathrm{Aut}(#1)}

\def\inv{{}^{-1}}
\def\sle{\subseteq}

\newtheorem{theorem}{Theorem}[section]
\newtheorem{proposition}[theorem]{Proposition}

\newtheorem{lemma}[theorem]{Lemma}
\newtheorem{remark}[theorem]{Remark}
\newtheorem{remarks}[theorem]{Remarks}

\renewcommand{\le}{\leqslant}
\renewcommand{\ge}{\geqslant}

\begin{document}

\title{Generating groups by conjugation-invariant sets}
\author{Valery Bardakov, Vladimir Tolstykh, Vladimir Vershinin}

\address{Valery Bardakov\\
Sobolev Institute of Mathematics\\
Siberian Branch of the Russian Academy of Sciences\\
630090 Novosibirsk\\
Russia}
\email{bardakov@math.nsc.ru}

\address{Vladimir Vershinin \\
Department of Mathematics \\
University Montpellier II \\
Place Eug\`ene Bataillon \\
34095 Montpellier Cedex 5 \\
France}
\email{vershini@math.univ-montp2.fr}
\address{Vladimir Tolstykh\\
Department of Mathematics\\
Yeditepe University\\
34755 Kay\i\c sda\u g\i \\
Istanbul\\
Turkey}
\email{vtolstykh@yeditepe.edu.tr}
\maketitle

\begin{abstract}
Let $S$ be a generating set of a group $G.$ We say that
$G$ has {\it finite width} relative to $S$
if $G=(S\cup S\inv)^k$ for a suitable natural number $k.$
We say that a group $G$ is a group of {\it finite $C$-width}
if $G$ has finite width with respect to all conjugation-invariant
generating sets of $G.$ We give a number of examples of
groups of finite $C$-width, and, in particular,
we prove that the commutator subgroup $F'$ of
Thompson's group $F$ is a group of finite $C$-width.
We also study the behaviour of the class of all
groups of finite $C$-width under some group-theoretic
constructions; it is established, for instance, that
this class is closed under formation of group
extensions.
\end{abstract}

\section{Introduction}

Let $G$ be a group and let $S$ be a generating
set of $G.$
Then the {\it length} $|g|_S$ of an element $g \in G$
with regard to $S$ is the least number of elements
of $S \cup S^{-1}$ whose product is $g$ and the  {\it width}
$\wid(G,S)$ of $G$ with regard to $S$
is
$$
\wid(G,S) = \sup \{|g|_S : g \in G\}.
$$
Thus $\wid(G,S)$ is either the least
natural number $k$ such that every element of $G$
is written as a product of at most
$k$ elements of $S \cup S^{-1},$
or $\wid(G,S)=\infty$ if such $k$ does not exist.
If the width $G$ with regard to $S$ is a
finite number $k,$ it is convenient to say
that $G$ {\it is generated by $S$ in $k$ steps.} If, furthermore,
$S=S^{-1}$ is a symmetric generating set,
then $G=S^k.$

Recall that a group $G$ is said to be a group {\it of finite
width} (or to satisfy the {\it Bergman property}),
if $G$ has finite width with respect to every
generating set.
The first example of an infinite group of finite
width has been given in 1980
by S.~Shelah \cite{Sh}. In 2003 Bergman (see \cite{Berg}) proved
that all infinite symmetric groups are groups of finite
width. Some questions Bergman included in his 2003 preprint of \cite{Berg}
initiated  the search of new examples of groups
of finite width. In particular, it has been discovered
that the following groups are groups
of finite width: the automorphism groups of 2-transitive
linearly ordered sets \cite{DrHo},
$\omega_1$-existentially closed groups \cite{Corn},
the autohomeomorphism groups of some topological spaces \cite{DrGo},
the automorphism groups of infinite-dimensional vector
spaces over arbitrary skew fields \cite{To_SMZh},
the automorphism groups of infinitely generated
free nilpotent groups \cite{To_JLM}, the automorphism
groups of many $\omega$-stable and  $\omega$-categorical
structures \cite{KeRo} and so on.

The class of all groups of finite width has a number of attractive
properties being, for instance, closed under homomorphic
images and under formation of group extensions \cite{Berg}.

Naturally, one can weaken the condition
requiring that a given group $G$ have finite width with regard to {\it all}
generating sets (which leads to the notion of a group
of finite width) by requiring that $G$ have finite
width with regard to all generating sets
satisfying a certain condition $\chi,$
thereby getting the notion of a group of finite
$\chi$-width.

In the present paper we study one of the
natural choices for $\chi$ above, the condition
$C=$ `to be conjugation-invariant' (to be invariant
under all inner automorphisms).
Thus we say that a given group
 $G$ is a {\it group of finite $C$-width}
if $G$ has finite width with regard to every
conjugation-invariant generating set.
Every group of finite width is of course a group
of finite $C$-width, but the converse is not
true (since, for instance, there are finitely
generated infinite groups of finite $C$-width).
The notion of a group of finite width
and the notion of a group of finite $C$-width
coincide on the class of all abelian
groups; Bergman \cite{Berg} proved that
the abelian groups of finite width are exactly
the finite ones.

It is easy to see that the special linear
group $\mathrm{SL}_n(K)$ where $n \ge 3$
and $K$ is an arbitrary field is a
group of finite $C$-width.
Indeed, the conjugacy class $T \sle \mathrm{SL}_n(K)$ which consists
of all transvections generates the group $G$
in $n$ steps: $\mathrm{SL}_n(F)=T^n.$
On the other hand, if $S=S^{-1}$ is a symmetric,
conjugation-invariant
generating set of $\mathrm{SL}_n(K),$
then an appropriate power $S^k$ contains
a transvection (and hence all transvections),
whence $\mathrm{SL}_n(K)=(S^k)^n=T^n=S^{kn}.$

In the next section we shall give some
other examples of groups of finite $C$-width.
In particular, we shall prove that
the commutator subgroup
$F'$ of Thompson's group $F$
is a group of finite $C$-width.
The proof uses the ideas developed
in the paper \cite{BIP} by Burago, Ivanov
and Polterovich. We shall also show
that some groups of `bounded' automorphisms
of infinitely generated relatively
free algebras are groups of finite $C$-width.

In the final section we analyze the behavior
of the class of all groups of finite $C$-width
under some group-theoretic constructions.
It shall be proved that the
class of all groups of finite $C$-width
is closed, like the class of all groups
of finite width, under formation
of group extensions. We shall consider
the question when a free product of nonidentity
groups is a group of finite $C$-width.
It turns out that the only, up to an isomorphism, free product
of nonidentity groups which is a group
of finite $C$-width is
the group $\Z_2 * \Z_2.$
We also give a number of necessary and sufficient
conditions for all functions $L : G \to \R$ on a given
group $G$ such that
\begin{alignat*} 3
& L^{-1}(0) &&=\{1\}, && \\
& L(g\inv) &&=L(g) &&\qquad (g\in G),\\
& L(gh) &&\le L(g)+L(h) && \qquad (g,h\in G), \\
& L(hgh\inv) &&= L(g) && \qquad (g,h\in G),
\end{alignat*}
called {\it norms} in \cite{BIP}, to be bounded
from above. For instance, all norms on a given group $G$
are bounded from above if and only if $G$ is a group
of finite $C$-width and every exhaustive chain
$(N_k)$ of normal subgroups of $G$ terminates
after finitely many steps. One more such a condition
states that for every action of $G$
by isometries on a metric space $\str{M;d}$ such that
$$
d(ghg\inv a,a) =d(ha,a) \qquad (a \in M; g,h \in G),
$$
the diameters of all $G$-orbits are bounded.

\section{Examples of groups of finite $C$-width}

We shall say that a generating set $S$
of a given group $G$ is {\it conjugation-invariant},
if $S$ is invariant under all conjugations
(inner automorphisms of $G$).
A group $G$ is said to be a group of
{\it finite $C$-width} if $G$ has finite
width relative to every conjugation-invariant
generating set.

The following simple sufficient condition
for finiteness of $C$-width works,
as we shall see, in many cases.

\begin{lemma} \label{trivlemma}
Let a group $G$ be generated
by finitely many conjugacy classes $C_1,\ldots,C_m$
in finitely many steps:
$$
\wid(G,C_1 \cup \ldots \cup C_m) = N < \infty.
$$
Then $G$ is a group of finite $C$-width.
\end{lemma}

\begin{proof}
Suppose that $C_k=a_k^G$ is the conjugacy
class of a certain element $a_k \in G$ $(k=1,\ldots,m).$
Consider a symmetric, conjugation-invariant generating set
$S$ of $G.$ For every $k=1,\ldots,m,$ a suitable power $S^{p_k}$ of $S$
contains the element $a_k.$ Set
$$
p=\max(p_1,\ldots,p_m).
$$
Then
$$
C_1 \cup \ldots \cup C_m \subseteq S^p,
$$
and it follows that $G=(S^p)^N=S^{pN}.$
\end{proof}

Now let us consider examples of groups of finite
$C$-width.

1) We have already mentioned in the introduction
that the special linear group $\mathrm{SL}_n(K)$
over a field $K$ of dimension $n \ge 3$ is a group
of finite $C$-width, since this group is generated
by the conjugacy class of all transvections
in $n$ steps.

2) Lemma \thelemma\ and Example 1.1 from
\cite{BIP} imply that the special linear
group  $\mathrm{SL}_n(\Z)$ where $n \ge 3$
over the ring of integers $\Z$ is also a group of finite $C$-width. In fact
this result has been proven in \cite{BIP} (though
in different terms). Indeed, it is known that the
group $\mathrm{SL}_n(\Z)$ is generated by the family
of all elementary transvections in some $K_n < \infty$ steps  (see \cite{Morr}).
An elementary transvection $t_{ij}(m),$ where $i,j$
with $1 \le i,j \le n$ are distinct indices
and $m \in \Z,$ is the matrix
$$
t_{ij}(m)=I+m I_{ij}.
$$
Observe that the group  $\mathrm{SL}_n(\Z)$ is also generated
by the conjugacy class $T$ of the transvection $t_{12}(1).$
On the other hand, the following well-known formula
$$
t_{ij}(m) = [t_{ik}(1), t_{kj}(m)] \qquad (m \in \Z),
$$
where $1 \le i,j,k \le n$ are pairwise distinct
indices, implies that every elementary transvection
is a product of two elements of $T.$ Thus
the group $\mathrm{SL}_n(\Z)$ is generated by
the conjugacy class $T$ in at most $2K_n$ steps.

3) Among known examples of groups of finite
width, quite a few are groups that are generated
by a {\it single} conjugacy class. For instance,
the mentioned property is shared by the infinite symmetric
groups, the automorphism groups of infinite
dimensional vector spaces over skew fields,
the automorphism groups of infinitely
generated free nilpotent groups and
so on.

4) Every nonidentity algebraically closed group $G$ is generated
by the conjugacy class of any nontrivial
element in two steps \cite[Cor. 2]{Macint}. Accordingly,
$G$ is a group of finite $C$-width.

5) Formally, if a given group has only
finitely many conjugacy classes, then
this group is a group of finite  $C$-width.

6) Evidently, any homomorphic image (any quotient
group) of a group of finite $C$-width
is also a group of finite $C$-width.

Let $G$ be a group. The commutator subgroup of a subgroup
$H$ of $G$ will be denoted by $H'.$ Our next example is

\subsection{The commutator subgroup $F'$ of Thompson's group $F$}
By the definition, subgroups $H_1,H_2$
of a group $G$ are said to be {\it commuting} if $[H_1,H_2]=\{1\},$
that is, if every element of $H_1$ commutes with every element
of $H_2.$

For convenience's sake, we reproduce some definitions
and results from \cite{BIP} we are going to use
below.

Let $H \le G,$ let $f$ be an element of
$G,$ and let $m \ge 1$ be a natural number. Then $f$ {\it
$m$-displaces} $H,$ if the subgroups
$$
H, fHf^{-1}, f^2 Hf^{-2},\ldots,f^{m} H f^{-m}
$$
are pairwise commuting \cite[Sect. 2.1]{BIP}.
We shall say that $f$ {\it $\omega$-displaces}
$H$ if $f$ $m$-displaces $H$ for all $m \ge 1.$

An element $x \in G$ is said to be an {\it $f$-commutator}
if $x$ is conjugate to a commutator
of the form $[f,g]$ where $g \in G.$ Observe
that the inverse of an $f$-commutator
is also an $f$-commutator.

\begin{proposition} \label{resultsfromBIP}
{\rm (i).} Suppose that $f$ $m$-displaces
a subgroup $H \le G$ for some $m \ge 2.$ Then
every element of $H'$ which has commutator
length $m$ in $H$ can be written
as a product of an $f$-commutator
and a commutator of some elements of $G$ {\rm
(see the proof of Theorem 2.2 (i)
in \cite{BIP})}.

{\rm (ii).} Suppose that $f \in G$ $\omega$-displaces
a subgroup $H$ of $G.$ The commutator length {\rm(}in $G${\rm)} of
every element $h \in H'$ is at most two {\rm \cite[Lemma 2.2]{BIP}}.

{\rm (iii).} Suppose that $f$ $m$-displaces
a subgroup $H$ of $G$ where $m \ge 2.$ Then
every commutator of elements of $H$ is a
product of two $f$-commutators {\rm \cite[Lemma 2.7]{BIP}}.
\end{proposition}

\begin{proposition} \label{F'_is_of_finite_c-width}
Let $f$ be an arbitrary nonidentity element of the
commutator subgroup $F'$ of Thompson's group
$F.$ Then every element of
$F'$ is a product of at most six elements
of $C(f)^{\pm 1}$ where $C(f)$ is the conjugacy
class of $f$ in $F'.$ Consequently,
the commutator subgroup $F'$ of Thompson's group $F$
is a group of finite $C$-width.
\end{proposition}

\begin{proof}
We shall work with the group $\mathrm{PL}_2([0,1]),$
one of the standard realizations of $F$ \cite{Belk,CFP}.
So we assume that $F=\mathrm{PL}_2([0,1]).$

The next result is well-known (folklore; see
Proposition 8.1 in \cite{AM}).

\begin{lemma}
Every commutator $[a,b] \in F'$ where $a,b \in F$ can
be written in the form
$[a_1,b_1]$ where $a_1,b_1$ are already elements
of $F'.$
\end{lemma}

Recall that every element of $F'$ acts
identically on a certain closed
segment $[0,\beta]$ and on a certain
closed segment
 $[\gamma,1]$ where $0 < \beta < \gamma <1$
(see, e.g., \cite[Cor. 1.3.4]{Belk}).
Let $f$ be a nonidentity element of
$F'.$ Suppose that $0 < \alpha_0 \in [0,1]$ is
the first nontrivial fixed point of $f,$ that is, the fixed point of $f$
such that the interval $[0,\alpha_0]$ is fixed
by $f$ pointwise, and every open neighbourhood
of $\alpha_0$ has points that are not fixed by $f.$
Let, further, $\alpha_1 < 1$ be the next fixed
point of $f,$ that is, the point of
$(0,1)$ such that there are no fixed points
of $f$ in the open interval $(\alpha_0,\alpha_1).$ Take a dyadic
point $\alpha \in (\alpha_0,\alpha_1) \cap \Z[\frac 12].$
Then either
$$
\alpha_0 < \alpha < f \alpha < f^2 \alpha < \ldots < f^n \alpha < \ldots
<
\alpha_1,
$$
in the case when $\alpha < f \alpha,$ or
$$
\alpha_0 < \ldots < f^n \alpha < \ldots < f^2 \alpha < f\alpha < \alpha
<
\alpha_1
$$
in the case when $\alpha > f \alpha,$ since $f$ is an increasing function.

Consider the open interval $J$ with endpoints
$\alpha$ and $f\alpha.$ Clearly, the intervals
\begin{equation}
J, fJ, f^2 J, \ldots, f^m J, \ldots
\end{equation}
all containing in the open interval $(\alpha_0,\alpha_1),$
are pairwise disjoint.

Let $H$ be the subgroup of all members $F'$
whose supports are in $J.$ It is easy to see
that $H$ is isomorphic to $F,$ since
$$
\mathrm{PL}_2(\overline J) \cong \mathrm{PL}_2([0,1])=F,
$$
where $\overline J$ is the closure of $J$
\cite[Prop. 1.4.4]{Belk}.

Recall that for every $n \ge 1,$ the group $F'$ acts
transitively on the family of $n$-element ordered
tuples of $(0,1) \cap \mathbf{Z}[1/2]$
\cite[Lemma 4.2]{CFP}. It follows that
every finite tuple of elements of
$F'$ can be taken into $H'$ by conjugation
by a suitable element of $F'.$

As open intervals in (\theequation)
are pairwise disjoint, we get that $f$ $\omega$-displaces $H$ in $F'.$
Using Lemma \thelemma\ and part (i) of Proposition \ref{resultsfromBIP},
we can write any element $h \in H'$
as a product of an $f$-commutator and
a commutator of elements of $F'$:
\begin{equation}
h = [f,b_1]^{c_1} [b_2,b_3].
\end{equation}
Conjugating then the elements $b_2,b_3$ by an appropriate
$c \in F',$ we obtain the commutator
$$
[b_2^c,b_3^c]
$$
of elements of $H.$ By
part (iii) of Proposition \ref{resultsfromBIP},
the latter commutator is a product of two $f$-commutators.
Then we deduce from (\theequation), that $h^c$
is a product of three $f$-commutators.
Accordingly, $h$ is a product of three
$f$-commutators, and then every element of
$F'$ is a product of three $f$-commutators.
Clearly, a product of three $f$-commutators
is a product of six elements of $C(f)^{\pm 1},$
where $C(f)$ is the conjugacy class of $f$ in $F'.$
\end{proof}

\begin{remarks}
\rm (i) As it has been demonstrated in the course
of the proof of the proposition, every
element of $F'$ is a product of two commutators
(folklore; see \cite[Prop. 8.1]{AM}).
The problem whether every element of $F'$
is a commutator is open.

(ii) The group $F$ itself is not a group of finite $C$-width,
for the abelianization $F/F' \cong \mathbf Z^2$
of $F$ is infinite. As for $F',$ it is not a group
of finite width, since its width relative to the
generating set
$$
\{x_n x_{n+1}^{-1} : n=0,1,2,\ldots\}
$$
is infinite.

(iii). Proposition \ref{F'_is_of_finite_c-width} implies the
 well-known fact of simplicity of the group $F'.$
\end{remarks}

\subsection{Subgroups of `bounded' automorphisms
of relatively free algebras} \label{aut(relfr)}

Let $\mathfrak V$ be a variety of algebras
and let $\cF \in \mathfrak V$ be a free $\mathfrak V$-algebra
of infinite rank $\vk.$ Consider a basis
$\cB$ of $\cF$ and an {infinite} cardinal $\l < \vk.$
For every $\f \in \mathrm{Aut}(\cF),$
define the {\it support} of $\f$ as
$$
\supp(\f) = \{b \in \cB : \f b \ne b\}.
$$
Set
$$
\mathrm{Aut}_{\l,\cB}(\cF) = \{\f \in \mathrm{Aut}(\cF) : |\supp(\f)|
\le
\l\}.
$$
Clearly, $G=\mathrm{Aut}_{\l,\cB}(\cF)$ is a subgroup of the automorphism group
$\mathrm{Aut}(\cF)$ of $\cF,$ and every element of $G$
fixes the `most' of elements of $\cB.$

\begin{proposition}
The group $G=\mathrm{Aut}_{\l,\cB}(\cF)$
is a group of finite $C$-width.
\end{proposition}

\begin{proof}
Recall that a subset $J$ of an infinite
set $I$ is called a  {\it moiety} of $I$ if
$|J|=|I \setminus J|.$

We term an automorphism of a given relatively
free algebra {\it permutational} if it
fixes setwise a certain basis of
this algebra.

Let $\cM$ be a relatively free algebra
with an infinite basis $\cX.$
Let $\cY$ be a moiety of $\cX.$ Consider
a permutational automorphism $\pi$ of
$\cM$ which acts on $\cY$ as an involution
without fixed points and which fixes
pointwise all elements of
$\cX \setminus \cY.$ As it has been demonstrated in the proof of part (i) of Theorem 1.5 in \cite{To_JLM},
every automorphism $\cM$
which fixes the subalgebra $\str\cY$ setwise
and takes every element of the set
$\cX \setminus \cY$ to itself is a product
of at most four conjugates of $\pi.$

Choose a subset $\cC$ of the basis $\cB$ of
cardinality $\l$ and a moiety $\cC_0$
of $\cC.$ Consider a permutational automorphism
$\pi^* \in G$ which
\begin{itemize}
\item fixes the set $\cB \setminus \cC_0$ pointwise;

\item acts on the set $\cC_0$ as an involution
without fixed points.
\end{itemize}

Let $H$ be a subgroup of $G$ consisting of
all automorphisms of $\cF$ preserving
the subalgebra $\str\cC$ generated by $\cC$ as a set
and taking to themselves all elements of $\cB \setminus \cC.$
Clearly, for every $\sigma \in G$
there is a permutational
automorphism $\rho \in G$ such that
the conjugate  $\sigma_0$ of $\sigma$ by $\rho$
is in $H,$ fixes the subalgebra $\str{\cC_0}$
setwise, and fixes the set $\cB \setminus \cC_0$
pointwise.

Then we obtain, as a corollary
of the above-quoted result from \cite{To_JLM},
that $\s_0$ is a product of at most
four conjugates of $\pi^*$ in $H \cong \aut{\str{\cC}}.$ Consequently, $\s$ is
a product of at most four conjugates
of $\pi^*$ in $G,$ whence $G=C(\pi^*)^4$ where $C(\pi^*)$
is the conjugacy class of $\pi^*$ in $G.$

Thus the group $G$ is generated
by the class $C(\pi^*)$ in at most
of four steps, and then it is a group
of finite $C$-width by Lemma \ref{trivlemma}.
\end{proof}

\section{$C$-width and group-theoretic constructions}

\subsection{Group extensions}
The proof of the fact that the class of all
groups of finite width is closed under
formation of group extensions given by Bergman in \cite{Berg}
is based upon the following statement.

\begin{lemma} \label{Bergman'sLemma}
\mbox{\rm \cite[Lemma 7]{Berg}}
Let $U=U^{-1}$ be a symmetric generating
set and let $H$ be a subgroup of given
group $G$ such that for some natural
number $n$
$$
Hg \cap U^n \ne \varnothing
$$
for all $g \in G.$ Then $H$ is generated by those
elements $G$ that have length $\le 2n+1$ relative to $U$:
$$
H =\str{H \cap U^{2n+1}}.
$$
\end{lemma}

\begin{proposition}
Let $G$ be a group and let $H$ be a normal
subgroup of $G$ such that both the subgroup $H$ and
the quotient group $G/H$ are groups
of finite $C$-width. Then $G$ is also
a group of finite $C$-width.
\end{proposition}

\begin{proof}
Let $U$ be a symmetric, conjugation-invariant
generating set of $G$ and let $f : G \to G/H$ be the
natural homomorphism from $G$ onto $G/H.$ Clearly, $f(U)$
is also a conjugation-invariant generating set of
the group $f(G)=G/H.$ But then
$$
\wid(f(G),f(U))= n < \infty,
$$
and we are in the conditions of Lemma \ref{Bergman'sLemma}.
Evidently, $H \cap U^{2n+1}$ is a conjugation-invariant
generating set of $H$:
$$
(H \cap U^{2n+1})^h = H^h \cap (U^{2n+1})^h =
H \cap U^{2n+1} \qquad (h \in H).
$$
Hence, by the conditions,
$$
\wid(H, H \cap U^{2n+1}) = m < \infty,
$$
whence
$$
\wid(G,U) \le \wid(f(G),f(U)) + \wid(H, U) \le
n + (2n+1)m < \infty.
$$
\end{proof}

We then obtain as a corollary that the class of all
groups of finite $C$-width is closed under
formation of group extensions, and, in particular,
it is closed under formation of direct products
and under formation of group extensions by finite
groups.

\subsection{Free products} Essentially, the following
result states that the only free product of nonidentity
groups which is a group of finite $C$-width is, up to an isomorphism,
the
free product $\Z_2 * \Z_2.$

\begin{proposition}
The free product $\prod^*_{i \in I} G_i$
of a family $\{G_i : i \in I\}$ of nonidentity
groups is a group of finite $C$-width
if and only if $|I|=2$ and both groups participating
in the free product are of order two.
\end{proposition}

\begin{proof}
$(\Leftarrow).$ Consider groups $A=\str a$ and $B =\str b,$
both isomorphic to the group $\Z_2.$ Let $w \in A*B$ be
a reduced word in letters $a,b$ which begins with $a.$
Then
$$
w=(ab)^k a, \text{ or } w=(ab)^k
$$
for a suitable natural number $k.$ It is easy
to see that each word of the form $(ab)^k a$ is conjugate
either to $a,$ or to $b.$ Let $k \ge 1.$ Then
$$
(ab)^k = (ab)^{k-1} ab = (ab)^{k-1} a \cdot b
$$
and the element $(ab)^k$ is a product of at most
two conjugates of elements of
$\{a,b\}.$ The argument in the case when $w$ begins
with $b$ is similar. Therefore the group
$A*B$ is generated in two steps
by the union of the conjugacy class of $a$ and
the conjugacy class of $b.$ Apply Lemma \ref{trivlemma}
to complete the proof.

$(\Rightarrow).$ Recall that a map $\Delta : H \to \Z$
from a given group $H$ into $\Z$ is called
a {\it quasi-homomorphism} if there is
a constant $C$ such that
$$
\Delta(ab) \le \Delta(a)+\Delta(b)+C \qquad (a,b \in H).
$$

Fix a family $\{G_i :i \in I\}$ of nonidentity groups,
and let $G$ denote the free product $\prod_{i \in I}^* G_i.$
An element $p$ of the free product $G$ is called a {\it palindrome}
if the reduced word representing $p$
(whose syllables are nonidentity
elements of factors $G_i$) is read the same way forwards and backwards. Thus
if $p \in G$ and
$$
p = v_1 \ldots v_n
$$
where $v_k$ are nonidentity
elements of free factors $G_i$
such that the elements $v_m,v_{m+1}$ lie in distinct
free factors for all $m,$ then $p$
is a palindrome if and only if
$$
v_1 \ldots v_n =v_n \ldots v_1.
$$
In \cite{BaTo} quasi-homomorphisms have been used
to show that all free products of nonidentity groups
that are not isomorphic to the group
$\Z_2 * \Z_2$ have infinite width
relative to the (generating) set of all palindromes,
or, in other words, infinite palindromic width.

It particular, it is proved in \cite{BaTo} that
if the free product $G=\prod_{i \in I}^* G_i$ where $|I| \ge 2$
is such that
\begin{equation}
\text{at least one free factors $G_i$ has an element of order $\ge 3$}
\end{equation}
then there is a quasi-homomorphism $\Delta_1 : G \to \Z$
with the following properties:
\begin{itemize}
\item[(a)] $\Delta_1(ab) \le \Delta_1(a)+\Delta_1(b) + 9$
for all $a,b \in G;$

\item[(b)] the value $\Delta_1$ at any palindrome of $G$
is at most $2;$

\item[(c)] $\Delta_1$ is not bounded from above.
\end{itemize}
(see \cite[pp. 203--204]{BaTo}).

1) Suppose then that our free product $G$
has the property (\theequation). The reader can verify
quite easily that it follows from the
definition of $\Delta_1$ in \cite{BaTo} that
$$
\Delta_1( a b a\inv) \le \Delta_1(b)+9 \qquad (a,b \in G).
$$
By (b), this implies that the value of $\Delta_1$
at any conjugate of a palindrome in $G$
is at most $11.$ Now, were the width of $G$ relative
to the family of all conjugates of palindromes
finite, the quasi-homomorphism $\Delta_1$ would be bounded
on $G,$ contradicting (c). Therefore every free product
of nonidentity groups with (\theequation) is not
a group of finite $C$-width.

2) Now let $G \not\cong \Z_2 * \Z_2$ not satisfy (\theequation). Then
$G$ is a free product of nontrivial
abelian groups of exponent two.
It is easy to see that in any free product
of abelian groups of exponent two, a conjugate
of a palindrome is a palindrome, too.
Thus the set of all palindromes of $G$ is invariant
under all conjugations. As we have mentioned
above, the palindromic width of $G$
is infinite, and hence $G$ is not a group
of finite $C$-width.
\end{proof}

\subsection{Functions similar to length functions
and cofinalities} Basing on the ideas
from the paper \cite{Berg} by Bergman, several authors
obtained a number of necessary and sufficient conditions
for all functions $L : G \to \R$
on a given group $G$ such that
\begin{alignat} 3
& L^{-1}(0) &&=\{1\}, && \\
& L(g\inv) &&=L(g) &&\qquad (g\in G), \nonumber \\
& L(gh) &&\le L(g)+L(h) && \qquad (g,h\in G). \nonumber
\end{alignat}
to be bounded from above. For instance, this takes
place if and only if
$G$ is a group of finite width and the cofinality $\mathrm{cf}(G)$
is uncountable \cite{DrHo, DrGo}; another such criterion states
that every action of $G$
by isometries on a metric space has bounded orbits \cite{Corn}.

By the definition, the {\it cofinality} $\mathrm{cf}(G)$
of an infinitely generated group $G$
is the least cardinal $\lambda$ such that $G$ can be written
as the union of a chain of cardinality $\lambda$ of its proper subgroups
(observe that no finitely generated group can be written
as the union of a chain of proper subgroups). In the case
when $\mathrm{cf}(G) > \aleph_0,$
a group $G$ is said to be a group of {\it uncountable
cofinality.} For instance, the symmetric
group $\sym X$ of an infinite set $X$
is a group of uncountable cofinality, since
$\mathrm{cf}(\sym X) > |X|$ \cite{MN, Berg}.

In accordance with the terminology introduced
in \cite{DrGo}, if $G$ is a group of finite
width and of uncountable cofinality, $G$ is
called a group of {\it strong uncountable cofinality.}

As we mentioned in the introduction, the paper \cite{BIP}
contains a number of results on functions (termed
{\it norms} in \cite{BIP}) $L : G \to \R$ on groups
satisfying the conditions (\theequation) and taking constant
values on conjugacy classes, that is,
satisfying the additional condition
$$
L(ghg\inv)=L(h) \qquad (g,h \in G).
$$

We shall provide below some necessary and sufficient
conditions for all norms on a given group $G$ to be bounded
from above. Our conditions are in fact naturally-weakened
versions of conditions, equivalent to strong uncountable
cofinality, that can be found in papers \cite{Corn,DrGo,Rosen}.

\begin{proposition} \label{strong_unc_norm_conf}
Let $G$ be a group. Then the following are equivalent:

{\rm (i).} $G$ is a group of finite $C$-width and
every exhaustive chain $(N_k)$
$$
N_0 \le N_1 \le \ldots \le N_k \le \ldots \le G
$$
of normal subgroups of $G$ {\rm(}every increasing chain of normal
subgroups whose union is $G${\rm)} terminates
after finitely many steps;

{\rm (ii).} Every exhaustive chain $(U_k)$
$$
U_0 \sle U_1 \sle \ldots \sle U_k \sle \ldots \sle G
$$
of subsets of $G$ such that for every $i \in \N$
\begin{itemize}
\item $U_i$ closed under taking inverses;
\item $U_i$ is conjugation-invariant;
\item the product $U_i U_i$ is contained in a suitable $U_k$
\end{itemize}
terminates after finitely many steps;

{\rm (iii).} Orbits of every action of $G$
by isometries on a metric space $\str{M,d}$
such that
$$
d(a, ghg\inv a) = d(a,ha) \qquad (a \in M; g,h \in G)
$$
have bounded diameters;

{\rm (iv).} Every function $L : G \to \R$
taking constant values on conjugacy classes
of $G$ and such that
\begin{itemize}
\item $L(g)=0$ if and only if $g=1;$
\item $L(gh) \le L(g)+L(h)$ for all $g, h \in G$
\end{itemize}
is bounded from above.
\end{proposition}

\begin{proof}
(i) $\Rightarrow$ (ii). Clearly, the chain of subgroups
of $G$ generated by sets $U_i,$
$$
\str{U_0} \le \str{U_1} \le \ldots \le \str{U_k} \le \ldots \le G,
$$
is an exhaustive chain of normal subgroups of $G.$ Then $G=\str{U_j}$
for a suitable natural number $j,$
and hence $U_j$ is a symmetric, conjugation-invariant
generating set of $G.$ As $G$ is a group of finite $C$-width,
$G=U_j^s.$ By the conditions on the chain $(U_k),$ the power
$U_j^s$ is contained in some $U_m$ for an appropriate
$m \in \N,$ whence $U_m=G.$

(ii) $\Rightarrow$ (iii). Let $a$ be an arbitrary
element of a metric space $M$ satisfying (iii).
Set
$$
U_n = \{g \in G : d(a,ga) \le n\} \qquad (n\in \N).
$$
It follows from (iii) that every
$U_n$ is conjugation-invariant.
Let $g,h \in U_n.$ Then we have that
$$
d(a,gh a) \le d(a,ga)+d(ga,gha) = d(a,ga)+d(a,ha) \le n+n=2n.
$$
Consequently, $U_n U_n \sle U_{2n}.$ As the
chain $(U_n)$ terminates, we get that
$G=U_m$ for some $m \in \N.$ Hence
$$
d(a,ga) \le m
$$
for all $g \in G.$ Thus the diameter of the
orbit $\{g a : g \in G\}$ of $a \in M$
is at most $2m.$

(iii) $\Rightarrow$ (iv). Let $a,b \in G.$ Set
$$
d(a,b)=L(ab\inv).
$$
It is easy to see that $d$ is a metric on $G$ satisfying
the conditions in (iii) for the left action $G$ on itself. Indeed,
we have that
\begin{align*}
& d(a,b)=0 \iff L(ab\inv)=0 \iff ab\inv=1 \iff a=b,\\
& d(a,b)=L(ab\inv)=L(ba\inv)=d(b,a),\\
& d(a,b)=L(ab\inv)=L(ac\inv \cdot cb\inv)\le
L(ac\inv)+L(cb\inv)=d(a,c)+d(c,b),\\
& d(ga,gb)=L(gab\inv g\inv)=L(ab\inv)=d(a,b),\\
& d(a,ghg\inv a)=L(ghg\inv)=L(h)=d(a,ha)
\end{align*}
for all $a,b,g,h \in G.$ Then the orbit of $1 \in G$ under the left
action of $G$ on itself has a bounded diameter $m \in \N,$ or
$$
L(g)= d(g1,1) \le m \qquad (g \in G).
$$

(iv) $\Rightarrow$ (i). Let $S=S^{-1}$ be a symmetric,
conjugation-invariant generating set of
$G.$ Then the function
$$
L_1(g) = |g|_S \qquad (g \in G)
$$
that is, the length function with regard to $S,$
which meets all conditions mentioned in (iv),
must be bounded from above by some natural
number $m.$ Accordingly, $G=S^m.$

Let further $(N_k)$ be an exhaustive chain
of normal subgroups of $G.$ For every $g \in G$ set
$$
L_2(g) =\min \{k \in \N : g \in N_k\} \qquad (g \in G).
$$
It is readily seen that $L_2$ satisfies all conditions in (iv). For example,
$$
L_2(ab) \le \max(L_2(a),L_2(b)) \le L_2(a)+L_2(b)
$$
for all $a,b \in G.$ One again concludes that $L_2$
is bounded from above by a certain natural number $m,$ whence $G=N_m.$
\end{proof}

\begin{remark}
\rm In the case when $G$ is a simple group, Proposition
\ref{strong_unc_norm_conf} provides a criterion of finiteness of $C$-width
of $G.$
\end{remark}

In the conclusion of the section we shall discuss some
notion which generalizes both the notion of a group
of finite width and the notion of a group of finite
$C$-width. Let $G$ be a group and let $\Sigma \le
\mathrm{Aut}(G)$
be a subgroup of the automorphism group of $G.$ We say
that $G$ has {\it finite $\Sigma$-width} if
$G$ has finite width with respect to
all $\Sigma$-invariant generating sets. Clearly, the case when $\Sigma=\{\mathrm{id}\}$
corresponds to the notion of a group of finite width, and the
case when $\Sigma=\mathrm{Inn}(G)$ to that one of finite $C$-width.

Intuitively, the greater the (setwise) stabilizer in $\aut G$
of a given set $S$ of generators of $G$, the
more `massive' $S$ appears to be with the `point of view'
of the automorphism group of $G.$
Thus if $G$ has finite $\Sigma$-width
in the case when $\Sigma=\mathrm{Aut}(G),$ it has
finite width with regard to all `most massive'
generating sets.

One can, as we did in Proposition \ref{strong_unc_norm_conf}, add
to the condition of finiteness of $\Sigma$-width
the condition of termination of all exhaustive chains of
$\Sigma$-invariant subgroups of $G.$ Modifying
then the formulation of Proposition \ref{strong_unc_norm_conf} accordingly,
one can obtain necessary and sufficient conditions for $G$ to have
finite $\Sigma$-width and, simultaneously, to satisfy the condition of termination
of all exhaustive chains of $\Sigma$-invariant
subgroups. For instance, the analogue of the
part (iii) of Proposition \ref{strong_unc_norm_conf} is
as follows: every action of $G$ by isometries
on a metric space $\str{M,d}$ for which
$$
d(\s(g)a,a)=d(ga,a), \qquad (a \in M, \s \in \Sigma, g \in G)
$$
has bounded orbits.

The case when $\Sigma$ is equal to the full
automorphism group $\mathrm{Aut}(G)$
of $G$ seems to be quite interesting.
Simplifying the terminology somewhat, we say
that $G$ has {\it finite {\rm Aut}-width} if $G$
has finite width with respect to every generating
set which is invariant under all automorphisms
of $G.$

Our final result shows that the class of all groups
of finite Aut-width does not have some attractive
properties that its counterparts, the classes
of all groups of finite width and all groups
of finite of $C$-width, have (in particular,
this class is not closed under homomorphic
images). Nevertheless, the property of having/not having
finite Aut-width can be used to distinguish
between the isomorphism types of groups.

\begin{proposition}
A free group $F$ is a group of finite
{\rm Aut}-width if and only if its rank
is infinite.
\end{proposition}

\begin{proof}
Suppose that $F$ is of finite rank. Consider
the (generating) set $P$ of all primitive elements of $F.$
Clearly, $P$ is invariant under all automorphisms
of $F,$ but the width of $F$ relative to
$P$ is infinite \cite[Th. 2.1]{BShT}.

Now let $F$ be of infinite rank. This time,
the width of $F$ with regard to the
set of all primitive elements is two \cite[Th. 2.1]{BShT}.
Consider a symmetric generating set $S$ of $F$
which is invariant under all automorphisms of
$F.$ Then a certain power $S^k$ of $S$ contains
a primitive element $p \in F.$ As $S$ is invariant
under automorphisms of $F,$ the said power of
$S$ contains all primitive elements. Consequently, $F=S^{2k}.$
\end{proof}

\end{document}